\documentclass[11pt,leqno]{amsart}
\usepackage[T1]{fontenc}
\usepackage{lmodern}
\usepackage{amsmath,amssymb,amsthm,mathtools}
\usepackage{graphicx}
\usepackage{tikz}
\usetikzlibrary{calc,intersections,arrows.meta}
\usepackage[margin=1in,centering]{geometry}
\usepackage{float}
\usepackage{flafter}
\usepackage[backend=bibtex,style=numeric,sorting=nyt,giveninits=true]{biblatex}
\usepackage[colorlinks=true,linkcolor=blue,citecolor=blue,urlcolor=blue]{hyperref}
\hypersetup{pdftitle={Self-dual plane curves and Appell's construction},
 pdfauthor={Timofey Fedorov},pdfsubject={Projective duality and contact transformations}}
\numberwithin{equation}{section}
\newtheorem{theorem}{Theorem}[section]
\newtheorem{proposition}[theorem]{Proposition}
\newtheorem{lemma}[theorem]{Lemma}

\newtheorem{remark}[theorem]{Remark}
\newtheorem*{theoremA}{Theorem A}
\newtheorem*{theoremB}{Theorem B}
\newtheorem{propositionC}{Proposition}

\theoremstyle{definition}
\newtheorem*{example*}{Example}
\newcommand{\Pp}{\mathbb P}
\newcommand{\A}{\mathbb A}
\newcommand{\C}{\mathbb C}
\newcommand{\Aut}{\operatorname{Aut}}
\newcommand{\id}{\operatorname{id}}

\newcommand{\Spec}{\operatorname{Spec}}
\newcommand{\dd}{\,d}

\title[Self-dual plane curves and Appell's construction]{Self-dual plane curves and Appell's construction}
\author{Timofey Fedorov}
\address{National Research University Higher School of Economics, Moscow, Russia.}
\email{tsfedorov@edu.hse.ru}

\begin{document}
\begin{abstract}
We generalize a construction introduced by Appell
in his 1894 paper and show that, over an algebraically
closed field of characteristic different from $2$,
every irreducible reflexive self-dual plane curve
can be obtained by an explicit construction starting
from a curve on a cyclic quotient of the affine plane.
We also show that, over an arbitrary field of
characteristic different from $2$, a smooth projective
geometrically connected curve is isomorphic to the
normalization of a geometrically integral reflexive
self-dual plane curve if and only if it admits a
nontrivial involution.

Over $\C$, we interpret Appell's original construction
in terms of a correspondence
\[
\Pp(T\Pp^2)\longleftarrow W^\circ
\longrightarrow\Pp(T\Pp^2),
\]
where $W^\circ$ is a smooth threefold and both maps
pull back the canonical contact distribution on
$\Pp(T\Pp^2)$ to the same contact distribution
on $W^\circ$.
\end{abstract}
\maketitle

\setcounter{tocdepth}{1}
\tableofcontents

\section*{Introduction}

\subsection*{Notation and conventions}
Unless otherwise specified, we work over an algebraically
closed field $k$ of characteristic different from $2$.
Let $V$ be a three-dimensional vector space over $k$,
and denote its dual by $V^\vee$.
We identify points of the dual projective plane
$\Pp(V^\vee)$ with the corresponding lines in $\Pp(V)$.
For a one-dimensional subspace $L$ of a vector space $W$,
we denote by $[L]$ the corresponding point of $\Pp(W)$.
For a nonzero vector $v\in W$, we put $[v]=[kv]$.

The incidence threefold
\[
\mathcal I=
\{(x,\ell)\in\Pp(V)\times\Pp(V^\vee):x\in\ell\}
\]
is canonically isomorphic to the projectivized tangent
bundle $\Pp(T\Pp(V))$ via the map
\[
(x,\ell)\longmapsto (x,[T_x\ell]).
\]
Let $\pi_1:\mathcal I\to\Pp(V)$ be the first projection.
The canonical contact structure on $\mathcal I$ is
given by the rank-two distribution $\mathcal H\subset T \mathcal I$
defined by
\[
\mathcal H_{(x,\ell)}
=(d\pi_1)_{(x,\ell)}^{-1}(T_x\ell)
\subset T_{(x,\ell)}\mathcal I.
\]

A morphism $g:M\to N$ between smooth varieties
of the same dimension endowed with contact
distributions $\mathcal H_M$ and $\mathcal H_N$
is called \emph{contact} if
\[
(dg)_p^{-1}((\mathcal H_N)_{g(p)})
=(\mathcal H_M)_p
\]
for every $p\in M$. We also call a rational map
contact if it is contact on a dense open subset
on which it is regular and the two contact
distributions are defined.

Let $M$ be a smooth threefold with contact
distribution $\mathcal H_M$. A curve
$\Lambda\subset M$ is called \emph{Legendrian} if
\[
T_p\Lambda\subset(\mathcal H_M)_p
\]
at every smooth point $p$ of $\Lambda$.

Let $C\subset\Pp(V)$ be an irreducible plane curve.
Its dual $C^\vee\subset\Pp(V^\vee)$ is the closure
of the image of the Gauss map
\[
\gamma_C:C\dashrightarrow\Pp(V^\vee),
\]
which sends a smooth point of $C$ to its tangent line.
The conormal curve $\widehat C\subset\mathcal I$
is defined by
\[
\widehat C=\overline{\{(x,\gamma_C(x)):
x\in C_{\mathrm{sm}}\}}.
\]

An irreducible plane curve $C$ of degree at least
two is called \emph{reflexive} if the exchange
of the two factors identifies $\widehat C$ with
$\widehat{C^\vee}$, where $V^{\vee\vee}$ is
canonically identified with $V$. In this case
\[
(C^\vee)^\vee=C,
\]
and the Gauss maps of $C$ and $C^\vee$ are
birational and inverse to each other.
By the Segre--Wallace criterion (see \cite[\S2]{hefez}),
$C$ is reflexive if and only if its Gauss map
is separable. In particular, in characteristic zero
every irreducible plane curve of degree at least
two is reflexive.

For an invertible linear map $A:V\to V^\vee$,
we denote by
\[
A_*:\Pp(V)\longrightarrow\Pp(V^\vee)
\]
the induced projective isomorphism.
We say that $C$ is \emph{self-dual} if
$C^\vee=A_*(C)$ for some such map $A$.

For a linear isomorphism $F:W_1\to W_2$ between
three-dimensional vector spaces, the induced projective
isomorphism $F_*$ determines a contact isomorphism
between the corresponding incidence threefolds,
given by
\[
\widehat F(x,\ell)
=(F_*(x),((F^{-1})^\vee)_*(\ell)).
\]

We call $\widehat F$ the \emph{conormal lift}
of $F_*$.

Let $B:V\times V\to k$ be a nondegenerate
symmetric bilinear form. The linear isomorphism
\[
b_B:V\longrightarrow V^\vee,\qquad
v\longmapsto B(v,\cdot),
\]
induces the polarity
\[
\rho_B:\Pp(V)\longrightarrow\Pp(V^\vee).
\]
The line $\rho_B(x)$ is called the polar line
of $x\in\Pp(V)$.
Let $q(v)=B(v,v)$ be the quadratic form associated
with $B$. We call the conic
\[
K=\{[L]\in\Pp(V):q|_L=0\}
\]
the \emph{isotropy conic} of $q$.

We say that $C$ is \emph{self-dual under the
polarity $\rho_B$} if $\rho_B(C)=C^\vee$.

\subsection*{Main results}

The study of self-dual plane curves goes back
to the nineteenth century. Appell \cite{appell} constructed
curves self-dual under a polarity using families
of conics self-dual under the same polarity,
noting the analogy with Moutard's construction
\cite{moutard} of surfaces invariant under inversion.

Suchar \cite{suchar} studied plane curves self-dual
under the polarity defined by a circle or its
composition with a rotation about the centre of
the circle.
Wear \cite{wear1919,wear1920} and Duncan \cite{duncan}
studied self-dual rational plane curves and their
symmetries.
Primrose \cite{primrose} showed that a general
plane quartic with two cusps and one node is
self-dual under two distinct polarities.
Oka \cite[Theorem~2.8]{oka} proved the uniqueness,
up to projective equivalence, of a self-dual plane
sextic with three $(3,4)$-cusps.

In the broader differential-geometric setting,
Arnold \cite[Problem~1994-17]{arnold} asked for
a description of curves in real projective spaces
of arbitrary dimension that are projectively
equivalent to their osculating duals.
Fuchs and Tabachnikov \cite[\S6]{fuchstabachnikov}
constructed self-dual curves in the real projective
plane using Legendrian curves satisfying monodromy
conditions.

The present paper has two main goals. The first is
to reformulate Appell's construction over $\C$ in
terms of a correspondence on the incidence threefold
$\Pp(T\Pp^2)$ given by two contact maps of degree two.
The second is to describe all irreducible reflexive
self-dual plane curves and to characterize the smooth
projective curves that occur as their normalizations.

For $h\in k^\times$, consider the automorphism
\[
R_h:\A^2_{u,s}\longrightarrow\A^2_{u,s},\qquad
(u,s)\longmapsto(hu,-h^2s).
\]
If $h$ is a root of unity of order $m$, consider
the rational map
\begin{equation}\label{eq:psi-h}
\psi_h:\A^2_{u,s}\dashrightarrow\A^2_{U,V},\qquad
\psi_h(u,s)=
\begin{cases}
\bigl(u^m,s^2/u^4\bigr),&\text{if $m$ is odd},\\[2mm]
\bigl(u^m,s/u^{m/2+2}\bigr),&\text{if $m$ is even}.
\end{cases}
\end{equation}
On the open set $u\ne0$, the map $\psi_h$ is
the quotient morphism by the cyclic group
$\langle R_h\rangle$ onto the open set $U\ne0$
in $\A^2_{U,V}$.

\begin{theoremA}
Let $C\subset\Pp^2$ be an irreducible plane curve
of degree at least two. The following conditions
are equivalent.
\begin{enumerate}
\item[\textup{(i)}] The curve $C$ is reflexive
and self-dual.

\item[\textup{(ii)}] There exist a root of unity
$h\in k^\times$ and an irreducible curve
$B\subset\A^2_{U,V}$ meeting the open set $U\ne0$
such that the reduced closure
$D\subset\A^2_{u,s}$ of
\[
\psi_h^{-1}(B\cap\{U\ne0\})
\]
is irreducible and $\dd u\ne0$ in
$\Omega_{k(D)/k}$, and the image of the morphism
$f:X\to\Pp^2$ given by
\[
f=
\bigl[1:u+\kappa:
s+u\kappa+\tfrac12(u^2-\kappa^2)\bigr],
\qquad
\kappa=-\frac12\frac{\dd s}{\dd u},
\]
is projectively equivalent to $C$. Here $X$ is
the normalization of the projective closure
of $D$.
\end{enumerate}

In situation \textup{(ii)}, the morphism $f$
is birational onto its image $C'=f(X)$, and
\[
(C')^\vee=(A_h)_*(C'),
\]
where $A_h$ is the linear isomorphism
represented by the matrix
\begin{equation}\label{eq:A-h}
A_h=\begin{pmatrix}
0&0&h^{-2}\\
0&-h^{-1}&0\\
1&0&0
\end{pmatrix}.
\end{equation}

Equivalently, one may start with an irreducible
curve $D\subset\A^2_{u,s}$ satisfying
$R_h(D)=D$ and $\dd u\ne0$ in $\Omega_{k(D)/k}$.
\end{theoremA}

\begin{theoremB}
Let $k$ be a field of characteristic different from $2$,
and let $X$ be a smooth projective geometrically
connected curve over $k$.
The following conditions are equivalent.
\begin{enumerate}
\item $X$ admits a nontrivial involution.
\item $X$ is isomorphic to the normalization of a
geometrically integral reflexive self-dual plane curve
of degree at least two.
\item $X$ is isomorphic to the normalization of a
geometrically integral reflexive plane curve of
degree at least two that is self-dual under a polarity.
\end{enumerate}
Moreover, for every nontrivial involution
$\sigma:X\to X$ there exist a morphism
$f:X\to\Pp^2$ birational onto its image and
a polarity $\rho_B$ such that
\[
f^\vee=\rho_B\circ f\circ\sigma,
\]
where $f^\vee$ is the composition of $f$ with the
Gauss map of $f(X)$, extended to $X$.
All the morphisms and isomorphisms in the statement
are defined over $k$.
\end{theoremB}

As an application of Theorem~A, we construct
self-dual plane curves of large genus relative to their degree.

\begin{propositionC}\label{prop:high-genus}
For every positive integer $b$ divisible by $8$,
there exists an irreducible complex plane curve
$C_b$ self-dual under the polarity $(A_1)_*$,
where $A_1$ is given by \eqref{eq:A-h}.
Its degree $d_b$ and geometric genus $g_b$ satisfy
\[
d_b=16b^3+4b^2-b,\qquad
g_b=2^{-10/3}d_b^{4/3}+O(d_b).
\]
\end{propositionC}

Section~1 presents Appell's and Myller's construction.
Section~2 proves Theorem~A and Proposition~C
and gives an example of a reflexive self-dual plane
curve that is not self-dual under any polarity.
Section~3 proves Theorem~B.

\medskip
\noindent\textbf{AI usage.}
The author used ChatGPT (OpenAI) to help develop and draft mathematical
arguments and examples, search the literature, edit the prose, and
prepare the \LaTeX{} source. The author is responsible for the final text,
including its mathematical claims, proofs, and references.

\medskip
\noindent\textbf{Acknowledgements.}
I am grateful to my advisor Serge Lvovski for his support and guidance.
The study has been funded within the framework of the HSE University Basic
Research Program (HSE-BR-2025-060).

\section{A contact-geometric interpretation of Appell's construction}

Throughout this section we work over $\C$.

\subsection{Appell's family of conics}
Fix a nondegenerate symmetric bilinear form $B$
on $V$ with associated quadratic form $q$,
polarity $\rho_B$, and isotropy conic $K$.

For a point $p\in\Pp(V)\setminus K$, choose a
nonzero vector $u\in V$ such that $[u]=p$.
Appell associates with $p$ the conic $\Sigma_p$
defined by the quadratic form
\begin{equation}\label{eq:appellconic}
F_u(x)=2B(u,x)^2-q(u)q(x),\qquad x\in V.
\end{equation}
Consider the $B$-orthogonal reflection
\[
r_p:V\longrightarrow V,\qquad
x\longmapsto x-\frac{2B(u,x)}{q(u)}u.
\]
Both $\Sigma_p$ and $r_p$ are independent of
the choice of $u$. Since
\[
F_u(x)=-q(u)B(r_p(x),x)
\]
and $r_p$ is an invertible linear map self-adjoint
with respect to $B$, the quadratic form $F_u$
is nondegenerate. Moreover, $(r_p)_*(\Sigma_p)
=\Sigma_p$.

For a nonzero vector $x\in V$ such that
$[x]\in\Sigma_p$, the differential of $F_u$
at $x$ is
\[
(\dd F_u)_x=-2q(u)B(r_p(x),\cdot).
\]
Thus the Gauss map of $\Sigma_p$ satisfies
\begin{equation}\label{eq:reflpole}
\gamma_{\Sigma_p}
=\rho_B\circ (r_p)_*|_{\Sigma_p}.
\end{equation}
Consequently, $\rho_B(\Sigma_p)=\Sigma_p^\vee$.

The conics $\Sigma_p$ and $K$ meet at the two
distinct points of $K\cap\rho_B(p)$.
Since $(r_p)_*$ fixes the polar line $\rho_B(p)$
pointwise and $\gamma_K=\rho_B|_K$, the conics
$\Sigma_p$ and $K$ have the same tangent lines
at these two points.

Let $T\subset\Pp(V)$ be an irreducible plane curve
different from $K$. At a smooth point $p\in T$,
choose a local parameter $t$ vanishing at $p$
and a regular local lift $u(t)$ of $T$ to
$V\setminus\{0\}$, so that $[u(0)]=p$.

Following Appell \cite{appell}, consider the equations
\begin{equation}\label{eq:characteristic}
F_{u(0)}(x)=0,\qquad
\left.\frac{\partial}{\partial t}
F_{u(t)}(x)\right|_{t=0}=0,
\end{equation}
where $x\in V$. Their common zero locus
$Z_p\subset\Pp(V)$ is independent of the choice
of the local parameter and the lift.

As shown in the next subsection, if $p\notin K$,
then $Z_p\setminus K$ consists of two distinct
points unless $T_pT$ is tangent to $K$,
in which case it is empty.

Appell's construction associates with $T$
the closed subset
\[
C_T=
\overline{\bigcup_{p\in T_{\mathrm{sm}}\setminus K}
(Z_p\setminus K)}
\subset\Pp(V).
\]

\begin{proposition}\label{prop:appell-classical}
If $C_T$ is an irreducible curve of degree
at least two, then
\[
\rho_B(C_T)=C_T^\vee.
\]
Conversely, for every irreducible plane curve
$C\ne K$ of degree at least two satisfying
$\rho_B(C)=C^\vee$, there exists an irreducible
plane curve $T\ne K$ such that $C=C_T$.
\end{proposition}

We prove both assertions in the next subsection.
If $C_T$ is reducible, its irreducible curve
components need not be self-dual under $\rho_B$.

\begin{figure}[!htbp]
\centering
\includegraphics[width=\textwidth]{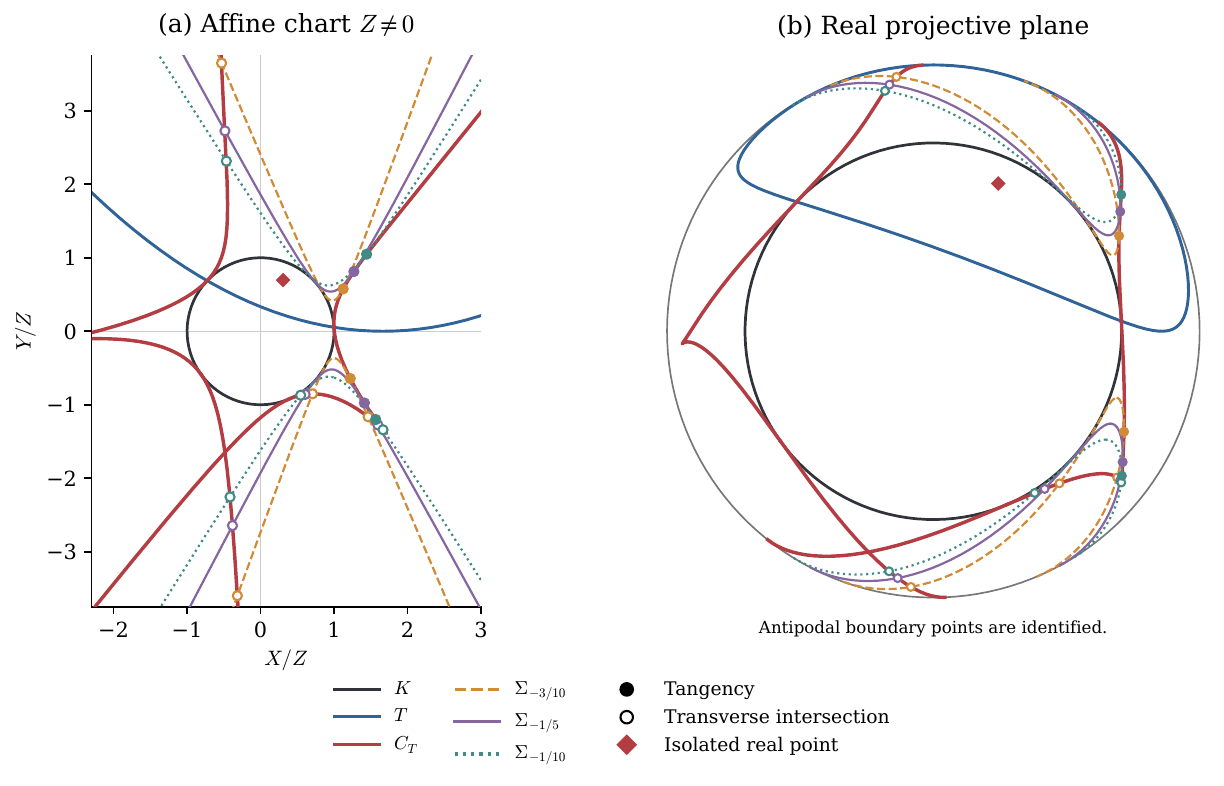}
\caption{Appell's construction for the parabola $T$.
The curves $K$, $T$ and $C_T$ are shown
in black, blue and red, respectively.
Filled circles indicate the points at which
the three selected conics are tangent to
$C_T$, and open circles indicate
their transverse intersections with
$C_T$.
The red diamond represents the isolated real
point of $C_T$.
In (b), the real projective plane is represented
by a disk with antipodal boundary points identified.}
\label{fig:appell-envelope}
\end{figure}

\begin{example*}
Let $q(X,Y,Z)=X^2+Y^2-Z^2$, and let $T$ be
the parabola defined by
\[
(3X-5Z)^2=75YZ.
\]
We parametrize $T$ by $[u(t)]$, where
\[
u(t)=\left(\frac53(1+t),\frac{t^2}{3},1\right).
\]
The pole of the tangent line to $T$ at $[u(t)]$
is $[v(t)]$, where
\[
v(t)=\left(-\frac{2t}{5},1,-\frac{t(t+2)}3\right).
\]
Put
\begin{align*}
a(t)&=q(u(t))
=\frac{t^4+25t^2+50t+16}{9},\\
b(t)&=q(v(t))
=-\frac{25t^4+100t^3+64t^2-225}{225}.
\end{align*}
The polynomials $a(t)$ and $b(t)$ have simple,
disjoint roots. Thus the double cover of $T$
defined by
\[
\lambda^2=\frac{a(t)}{b(t)}
\]
is irreducible. By
Proposition~\ref{prop:appell-contact},
the curve $C_T$ is irreducible and satisfies
\[
\rho_B(C_T)=C_T^\vee.
\]

The equation of $C_T$ can be obtained
by eliminating $t$ from \eqref{eq:characteristic}.
Writing $z=(X,Y,Z)\in V$ and regarding $F_{u(t)}(z)$
as a binary quartic by homogenizing with respect to $t$,
we compute
\[
\operatorname{Disc}_t F_{u(t)}(z)
=q(z)^2 H(z),
\]
where $H$ is a homogeneous polynomial of degree
eight, not divisible by $q$.

If $q(z)\ne0$, a multiple root of the binary
quartic $F_{u(t)}(z)$ cannot be a zero of $a(t)$,
since the roots of $a(t)$ are simple.
Hence the zero locus of $H$ in
$\Pp(V)\setminus K$ is contained in
$C_T$.
Since $C_T$ is irreducible, $H$ is
a constant multiple of a power of its defining
polynomial. A direct computation gives
\[
H\left(\frac53,\frac43,1\right)=0,\qquad
\frac{\partial H}{\partial X}
\left(\frac53,\frac43,1\right)\ne0.
\]
Consequently, $C_T$ is an irreducible
curve of degree eight.

Let $\alpha_1<\alpha_2$ be the real roots of $a$
and $\beta_1<\beta_2$ the real roots of $b$.
They satisfy
\[
\beta_1<\alpha_1<\alpha_2<\beta_2.
\]
For $t$ in either of the intervals
$(\beta_1,\alpha_1)$ and $(\alpha_2,\beta_2)$,
the two points
\[
p_\pm(t)=
\left[
\sqrt{b(t)}\,u(t)\pm\sqrt{a(t)}\,v(t)
\right]
\]
belong to $C_T$.
The two parametrizations extend continuously
to the endpoints of each interval, where their
values coincide.

Figure~\ref{fig:appell-envelope} shows the conics
$\Sigma_{[u(w)]}$ for
\[
w\in\left\{-\frac3{10},-\frac15,-\frac1{10}\right\}.
\]
Each of these conics is tangent to
$C_T$ at $p_+(w)$ and $p_-(w)$
and intersects it transversely at four other
real points.

To verify these assertions, put
\[
A=B(u(w),u(t)),\qquad
D=B(u(w),v(t)),\qquad c=a(w).
\]
Substituting the vectors
$\sqrt{b(t)}\,u(t)\pm\sqrt{a(t)}\,v(t)$
into $F_{u(w)}$ and eliminating the square roots
gives
\[
(bA^2+aD^2-cab)^2-4abA^2D^2=0.
\]
For each of the three values of $w$, the
left-hand side is divisible by $(t-w)^4$.
The quotient is a square-free polynomial of
degree twelve, nonzero at $t=w$, with four real
roots, all lying in
$(\beta_1,\alpha_1)\cup(\alpha_2,\beta_2)$.
A direct computation also shows that $AD\ne0$
at these four roots and that $\nabla H\ne0$ at
the corresponding points of $C_T$ and at
$p_+(w)$ and $p_-(w)$.
The four simple roots give four transverse
intersections, and the factor $(t-w)^4$ shows
that the intersection multiplicity at each of
$p_+(w)$ and $p_-(w)$ is exactly two.
The polynomials $a$ and $A$ are coprime,
as are $b$ and $D$, so no intersection lies on $K$.

The curve $C_T$ also has an isolated
real point with approximate affine coordinates
\[
(0.30643437,\;0.69687864).
\]
This point is marked by a diamond in
Figure~\ref{fig:appell-envelope}.
\end{example*}

\subsection{A correspondence on the incidence threefold}

The polarity $\rho_B$ determines an involution
\[
\tau:\mathcal I\longrightarrow\mathcal I,\qquad
(x,\ell)\longmapsto
(\rho_B^{-1}(\ell),\rho_B(x)).
\]
Locally on $\mathcal I$, choose regular vector-valued
lifts $u,v$ of $x$ and $\rho_B^{-1}(\ell)$,
respectively. Then
\[
B(u,v)=0,
\]
and the contact distribution $\mathcal H$ is
defined by the $1$-form $B(v,\dd u)$.
Since $B(v,\dd u)=-B(u,\dd v)$, the involution
$\tau$ preserves $\mathcal H$.

Let $W\subset\Pp(V\oplus V)$ be the subvariety
defined by the homogeneous quadratic equations
\[
B(u,v)=0,\qquad q(u)=q(v),
\]
where $(u,v)\in V\oplus V$ are the homogeneous
vector coordinates. We denote the point represented
by a nonzero pair $(u,v)$ by $[u:v]$.
Consider the rational map
\[
\pi:W\dashrightarrow\mathcal I,\qquad
[u:v]\longmapsto ([u],\rho_B([v])).
\]

Let
\[
\mathcal I^\circ=
\{(x,\ell)\in\mathcal I:
x\notin K,\ \rho_B^{-1}(\ell)\notin K\},
\]
and let $W^\circ\subset W$ be the complement
of the hypersurface defined by $q(u)=0$.
The restriction
\[
\pi:W^\circ\longrightarrow\mathcal I^\circ
\]
is an étale double cover. Indeed, for local lifts
$u,v$ of $x$ and $\rho_B^{-1}(\ell)$, respectively,
the preimages of $(x,\ell)\in\mathcal I^\circ$
are represented by pairs $(u,\lambda v)$ with
\[
\lambda^2=\frac{q(u)}{q(v)}.
\]
The ratio $q(u)/q(v)$, regarded as a rational
function on $\mathcal I$ after choosing rational
lifts, is not a square: it has a simple zero
along the divisor $\pi_1^{-1}(K)$ for lifts
regular and nonvanishing at its general point.
Thus $W^\circ$ is irreducible.

Since two orthogonal isotropic vectors in $V$
are linearly dependent, every point of
$W\setminus W^\circ$ is represented by a pair
$(az,bz)$, where $z\in V\setminus\{0\}$ satisfies
$q(z)=0$ and $(a,b)\ne(0,0)$.
Hence $W\setminus W^\circ$ is the image of
$\Pp^1\times K$ under the Segre embedding
\[
([a:b],[z])\longmapsto[az:bz].
\]
In particular, it has dimension two. Every
irreducible component of $W$ has dimension at
least three, since $W$ is defined by two equations
in $\Pp(V\oplus V)$. Consequently, $W$ is
irreducible.

The linear maps on $V\oplus V$ given by
\[
(u,v)\longmapsto(u,-v),\qquad
(u,v)\longmapsto(v,u),\qquad
(u,v)\longmapsto(u+v,u-v)
\]
induce projective involutions $\delta$, $\iota$
and $H$, respectively, preserving $W$ and
$W^\circ$. The involution $\delta$ is the deck
involution of $\pi:W^\circ\to\mathcal I^\circ$.
We have
\begin{equation}\label{eq:keyH}
\begin{gathered}
H\circ\delta\circ H=\iota,\qquad
\pi\circ\iota=\tau\circ\pi,\\
(\pi\circ H)\circ\delta=\tau\circ(\pi\circ H).
\end{gathered}
\end{equation}

Since $\pi:W^\circ\to\mathcal I^\circ$ is étale,
the pullback $\pi^*\mathcal H$ is a contact
distribution on $W^\circ$.
The involution $H$ preserves this distribution.
Indeed, for a local lift $(u,v)$ of a point
of $W^\circ$, differentiating the equations
$q(u)=q(v)$ and $B(u,v)=0$ gives
\[
B(u-v,\dd(u+v))=-2B(v,\dd u).
\]

Consider the rational map
\[
\pi\circ H:W\dashrightarrow\mathcal I,\qquad
[u:v]\longmapsto([u+v],\rho_B([u-v])).
\]
Both $\pi$ and $\pi\circ H$ restrict to étale
double covers of $\mathcal I^\circ$ by $W^\circ$
and are contact with respect to $\pi^*\mathcal H$
and $\mathcal H$. Their deck involutions are
$\delta$ and $\iota$, respectively.
We thus obtain the contact correspondence
\begin{equation}\label{eq:appell-correspondence}
\mathcal I^\circ
\xleftarrow{\ \pi\ }W^\circ
\xrightarrow{\ \pi\circ H\ }\mathcal I^\circ.
\end{equation}

An irreducible closed Legendrian curve
$\Lambda\subset\mathcal I$ with nonconstant first
projection is the conormal curve of its image
$C=\pi_1(\Lambda)$.
Indeed, at a general point $(x,\ell)\in\Lambda$,
the Legendrian condition gives $T_xC=T_x\ell$,
so $\ell$ is the tangent line to $C$ at $x$.

\begin{proposition}\label{prop:appell-contact}
Let $T\subset\Pp(V)$ be an irreducible plane curve
different from $K$. Put
$\Lambda=\widehat T\cap\mathcal I^\circ$ and
\[
\mathcal L=
\overline{(\pi\circ H)(\pi^{-1}(\Lambda))}
\subset\mathcal I,
\]
where the inverse image is taken in $W^\circ$.
Then $\tau(\mathcal L)=\mathcal L$, every
irreducible component of $\mathcal L$ is
Legendrian, and
\[
C_T=\pi_1(\mathcal L).
\]
If $\pi^{-1}(\Lambda)$ is nonempty and
irreducible, then $C_T$ is an
irreducible plane curve of degree at least two
satisfying
\[
\rho_B(C_T)=C_T^\vee.
\]
\end{proposition}

\begin{proof}
Let $p\in T\setminus K$ be a smooth point.
Choose a local parameter $t$ vanishing at $p$
and a holomorphic lift $u(t)$ such that
$q(u(t))=1$. Put $u=u(0)$ and $u'=u'(0)$.
The second equation of \eqref{eq:characteristic}
takes the form
\[
4B(u,x)B(u',x)=0.
\]
The hyperplane defined by $B(u,x)=0$ intersects
$\Sigma_p$ at the two points of $\Sigma_p\cap K$.

Choose a nonzero vector $v\in V$ such that
$\rho_B([v])=T_pT$. Then
\[
B(u,v)=B(u',v)=B(u,u')=0.
\]
Since $u$ and $u'$ are linearly independent,
the hyperplane defined by $B(u',x)=0$ is
$\langle u,v\rangle$.

If $q(v)=0$, the restriction of $F_u$ to
$\langle u,v\rangle$ vanishes only along
the subspace spanned by $v$. Thus the common
zero locus of \eqref{eq:characteristic} is
contained in $K$.

If $q(v)\ne0$, rescale $v$ so that $q(v)=1$.
Substituting $x=au+bv$ into
\eqref{eq:appellconic}, we obtain
\[
F_u(au+bv)=a^2-b^2.
\]
Hence the common zero locus of
\eqref{eq:characteristic} in $\Pp(V)\setminus K$
is $\{[u+v],[u-v]\}$.
By \eqref{eq:reflpole}, the tangent lines to
$\Sigma_p$ at these points are $\rho_B([u-v])$
and $\rho_B([u+v])$, respectively.

On the other hand,
\[
\pi^{-1}(p,T_pT)=\{[u:v],[u:-v]\},
\]
and
\begin{align*}
(\pi\circ H)([u:v])
&=([u+v],\rho_B([u-v])),\\
(\pi\circ H)([u:-v])
&=([u-v],\rho_B([u+v])).
\end{align*}
Since $\pi:W^\circ\to\mathcal I^\circ$ is
an étale double cover, the preimages of the
points of $\Lambda$ corresponding to smooth
points of $T$ are dense in $\pi^{-1}(\Lambda)$.
As $\pi_1$ is proper, we obtain
$C_T=\pi_1(\mathcal L)$.

The curve $\widehat T$ is Legendrian.
Since both maps in \eqref{eq:appell-correspondence}
are étale and contact, every irreducible component
of $\mathcal L$ is Legendrian.
Since $\delta$ preserves $\pi^{-1}(\Lambda)$,
\eqref{eq:keyH} implies that $\tau(\mathcal L)=\mathcal L$.

If $\pi^{-1}(\Lambda)$ is nonempty and
irreducible, then $\mathcal L$ is an
irreducible curve. Neither of its projections
can be constant, since otherwise
$\tau(\mathcal L)=\mathcal L$ would imply
that both projections are constant.
Consequently,
$\mathcal L=\widehat{C_T}$,
and $C_T$ is not a line.
The equality
$\tau(\widehat{C_T})
=\widehat{C_T}$
implies
$\rho_B(C_T)=C_T^\vee$.
\end{proof}

\begin{proof}[Proof of Proposition~\ref{prop:appell-classical}]
Suppose first that $C_T$ is irreducible of degree
at least two. Let $\mathcal L$ be as in
Proposition~\ref{prop:appell-contact}.
Then $\widehat{C_T}$ is an irreducible
component of $\mathcal L$. Since
$\tau(\mathcal L)=\mathcal L$, we have
\[
\rho_B^{-1}(C_T^\vee)
=\pi_1(\tau(\widehat{C_T}))
\subset C_T.
\]
Both sides are irreducible curves, so
$\rho_B(C_T)=C_T^\vee$.

Conversely, let $C\ne K$ be an irreducible plane curve
of degree at least two satisfying $\rho_B(C)=C^\vee$.
Then $\tau(\widehat C)=\widehat C$, and the open subset
$\Lambda=\widehat C\cap\mathcal I^\circ$ is nonempty.

Consider the closed subset
\[
Z=(\pi\circ H)^{-1}(\Lambda)\subset W^\circ.
\]
Since $\iota$ is the deck involution of
$\pi\circ H$, we have $\iota(Z)=Z$.
By \eqref{eq:keyH}, it follows that
$\tau(\pi(Z))=\pi(Z)$.

Let $\mathcal M\subset\mathcal I$ be the closure
of $\pi(Z)$. Since both maps in
\eqref{eq:appell-correspondence} are étale
and contact, every irreducible component of
$\mathcal M$ is a Legendrian curve.
Moreover, $\tau(\mathcal M)=\mathcal M$.
Hence at least one irreducible component of
$\mathcal M$ has nonconstant first projection.
Such a component is the conormal curve
$\widehat T$ of an irreducible plane curve
$T\ne K$.

It remains to show that $C=C_T$.
Since $\tau(\Lambda)=\Lambda$, \eqref{eq:keyH}
gives $\delta(Z)=Z$. As $\delta$ is the deck
involution of $\pi$, this gives
\[
\pi^{-1}(\pi(Z))=Z.
\]
Since $\pi$ is finite on $W^\circ$, the subset
$\pi(Z)$ is closed in $\mathcal I^\circ$.
Consequently,
\[
\widehat T\cap\mathcal I^\circ
\subset\pi(Z),
\]
and therefore
\[
(\pi\circ H)
\bigl(\pi^{-1}(\widehat T\cap\mathcal I^\circ)\bigr)
\subset\Lambda.
\]
The left-hand side is a nonempty closed curve
in $\mathcal I^\circ$, since both maps in
\eqref{eq:appell-correspondence} are finite
and surjective. By irreducibility of $\Lambda$,
the inclusion is an equality.
Taking closures in $\mathcal I$ and applying
Proposition~\ref{prop:appell-contact}, we obtain
$C_T=C$.
\end{proof}

\subsection{Myller's construction}

Myller \cite{myller} gave a geometric construction
of $C_T$.
Let $p\in T\setminus K$ be a smooth point such
that the tangent line $\ell=T_pT$ is not tangent
to $K$. Let $E$ and $F$ be the points at which
the two tangents from $p$ to $K$ meet $K$, and let
$\{R,S\}=\ell\cap K$.
The common zero locus of \eqref{eq:characteristic}
in $\Pp(V)\setminus K$ consists of the points
\[
P=ER\cap FS,\qquad Q=ES\cap FR.
\]
The chord $PQ$ passes through $p$.
If $L=\ell\cap EF$, then the tangent lines
to $\Sigma_p$ at $P$ and $Q$ are
\[
LP=\rho_B(Q),\qquad LQ=\rho_B(P),
\]
respectively.

\begin{figure}[!htbp]
\centering
\includegraphics[width=0.88\textwidth]{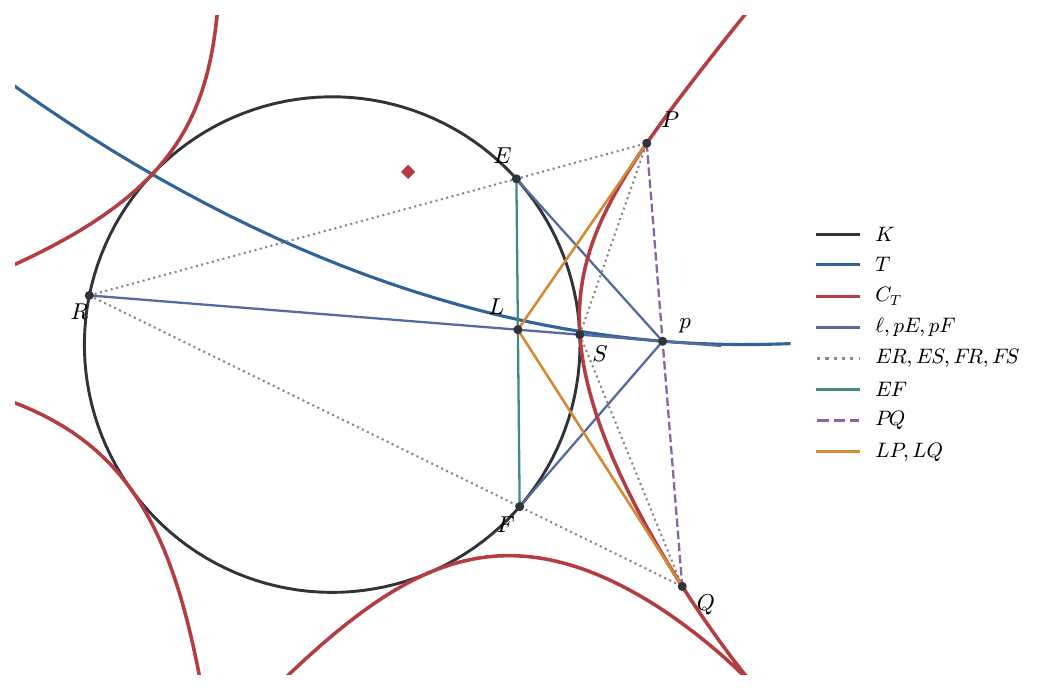}
\caption{Myller's construction for the parabola $T$
in Figure~\ref{fig:appell-envelope}, at the point
$p=[100:1:75]$.
The conic $\Sigma_p$ is tangent to
$C_T$ at $P$ and $Q$.
The lines $LP$ and $LQ$ are tangent to
$C_T$ at these points, respectively.}
\label{fig:myller}
\end{figure}

The lines $\ell$ and $PQ$ are exchanged by the
nontrivial involution of the pencil of lines
through $p$ fixing $pE$ and $pF$. Equivalently,
the cross-ratio of the four lines satisfies
\[
(pE,pF;\ell,PQ)=-1
\]
(see Figure~\ref{fig:myller}).

To verify these assertions, choose homogeneous
coordinates in which
\[
K=\{XZ-Y^2=0\},\qquad
E=[1:0:0],\qquad F=[0:0:1],\qquad p=[0:1:0].
\]
Since $\ell$ is not tangent to $K$, its equation
has the form $X=a^2Z$ for some $a\ne0$.
After ordering the points $R$ and $S$, we have
\[
R=[a^2:a:1],\qquad S=[a^2:-a:1].
\]
Consequently,
\[
P=[-a^2:a:1],\qquad Q=[-a^2:-a:1].
\]
The vectors
\[
u=(0,a,0),\qquad v=(-a^2,0,1)
\]
satisfy $[u]=p$, $\rho_B([v])=\ell$,
$B(u,v)=0$ and $q(u)=q(v)\ne0$.
Since $P=[u+v]$ and $Q=[u-v]$, both points
satisfy \eqref{eq:characteristic} by the proof
of Proposition~\ref{prop:appell-contact}.

Moreover, $L=[a^2:0:1]$, and the polar lines
$\rho_B(Q)$ and $\rho_B(P)$ have equations
\[
X+2aY-a^2Z=0,\qquad X-2aY-a^2Z=0,
\]
respectively. Hence $\rho_B(Q)=LP$ and
$\rho_B(P)=LQ$.

Finally, parametrize the pencil of lines through
$p$ by $[\lambda:\mu]\in\Pp^1$, assigning to
$[\lambda:\mu]$ the line $\mu X=\lambda Z$.
The lines $pE$, $pF$, $\ell$ and $PQ$ correspond
to the points
\[
[1:0],\qquad [0:1],\qquad
[a^2:1],\qquad [-a^2:1],
\]
respectively. Their cross-ratio is $-1$.

\subsection{The quotient by the involution
\texorpdfstring{$\tau$}{tau}}

The involution $\tau$ acts freely on
$\mathcal I^\circ$. Indeed, if
$\tau(x,\ell)=(x,\ell)$, then $\ell=\rho_B(x)$
and the incidence relation $x\in\ell$ implies
$x\in K$. Let
\[
J=\mathcal I^\circ/\langle\tau\rangle,\qquad
\eta:\mathcal I^\circ\longrightarrow J
\]
be the quotient and the quotient morphism.
Since $\tau$ preserves the contact distribution
on $\mathcal I^\circ$, there is a contact
distribution on $J$ with respect to which
$\eta$ is contact.

The involutions $\delta$ and $\iota$ generate
a group of order four acting freely on $W^\circ$.
By \eqref{eq:keyH}, the composition
\[
\eta\circ\pi:W^\circ\longrightarrow J
\]
is the quotient morphism by this group.
Since $H$ conjugates $\delta$ to $\iota$
and $\iota$ to $\delta$, it induces an
involution
\[
\overline H:J\longrightarrow J
\]
satisfying
\[
\overline H\circ\eta\circ\pi
=\eta\circ\pi\circ H.
\]
Since $H$ and $\eta\circ\pi$ are contact
and $\eta\circ\pi$ is étale, the involution
$\overline H$ is contact.

Choose a basis of $V$ in which the bilinear
form $B$ is represented by the matrix $A_1$
in \eqref{eq:A-h}.
Consider the affine chart of $\mathcal I$
consisting of pairs
\begin{equation}\label{eq:incidence-chart}
\bigl([1:x:z],[xv-z:-v:1]\bigr),
\qquad (x,z,v)\in\A^3.
\end{equation}
The contact distribution on this chart is
defined by the $1$-form $\dd z-v\,\dd x$.
Put
\begin{equation}\label{eq:affine-coordinates}
u=\frac{x+v}{2},\qquad
\kappa=\frac{x-v}{2},\qquad
s=z-u\kappa-\frac{u^2-\kappa^2}{2}.
\end{equation}
Then
\begin{equation}\label{eq:affine-contact-form}
\dd z-v\,\dd x=\dd s+2\kappa\,\dd u,
\end{equation}
and $\tau$ acts by
\[
\tau(u,s,\kappa)=(u,-s,-\kappa).
\]
Put $w=s^2$ and $p=-4s\kappa$. Then
\[
\C[u,s,s^{-1},\kappa]^{\langle\tau\rangle}
=\C[u,w,w^{-1},p].
\]
The intersection of the open set $s\ne0$
in the affine chart above with $\mathcal I^\circ$
is preserved by $\tau$. Its quotient is an open
subset of $\A^3_{u,w,p}$ on which the quotient
morphism is given by
\begin{equation}\label{eq:eta}
\eta(u,s,\kappa)=(u,s^2,-4s\kappa).
\end{equation}
Moreover,
\[
\dd w-p\,\dd u=2s(\dd s+2\kappa\,\dd u).
\]
Thus $J$ is rational, and the induced contact
distribution on $J$ is defined in these coordinates
by the $1$-form $\dd w-p\,\dd u$.

\begin{proposition}\label{prop:quotient}
The morphism
\[
(\pi,\pi\circ H):
W^\circ\longrightarrow
\mathcal I^\circ\times\mathcal I^\circ
\]
induces an isomorphism of $W^\circ$ with the
fibre product of the two morphisms
\[
\overline H\circ\eta:\mathcal I^\circ\to J,
\qquad
\eta:\mathcal I^\circ\to J.
\]
In particular, for a Legendrian curve
$\Lambda\subset\mathcal I^\circ$, we have
\[
(\pi\circ H)(\pi^{-1}(\Lambda))
=\eta^{-1}\bigl(\overline H(\eta(\Lambda))\bigr).
\]
\end{proposition}

\begin{proof}
By the definition of $\overline H$, the morphism
$(\pi,\pi\circ H)$ factors through the fibre product.
Both $W^\circ$ and the fibre product are finite
étale covers of $J$ of degree four.

Suppose that two distinct points of $W^\circ$
have the same image in the fibre product.
Since their images under $\pi$ coincide,
one is obtained from the other by applying $\delta$.
By \eqref{eq:keyH}, their images under $\pi\circ H$
are exchanged by $\tau$ and hence distinct, since
$\tau$ acts freely on $\mathcal I^\circ$.
This is a contradiction. Hence the morphism to the
fibre product is bijective on every fibre
over $J$ and is therefore an isomorphism.
The equality for $\Lambda$ follows from the
definition of the fibre product.
\end{proof}

Let us explain the relation between Appell's
construction and Theorem~A. The action of $\tau$ on
$(u,s,\kappa)$ is the case $h=1$ of \eqref{eq:diag}.

Let $C$ be an irreducible plane curve of degree
at least three satisfying $\rho_B(C)=C^\vee$.
Choose $T$ as in the proof of
Proposition~\ref{prop:appell-classical} and put
\[
M=\overline H(\eta(\widehat T\cap\mathcal I^\circ))
\subset J.
\]
Then
\[
\eta^{-1}(M)=\widehat C\cap\mathcal I^\circ.
\]
The function $u$ is nonconstant on $\widehat C$:
otherwise, the Legendrian condition and the
invariance under $\tau$ would imply that $s=0$,
and $C$ would be a conic. Moreover, $s$ does
not vanish identically on $\widehat C$, since
otherwise $\kappa=0$ and $C=K$.

Let $B_T\subset\A^2_{U,V}$ be the closure of the
image of $M$ under the rational map
\[
(u,w,p)\longmapsto(u,w/u^4).
\]
Since $M$ is Legendrian, we have
$p=\dd w/\dd u$, so this map is birational
onto $B_T$. By \eqref{eq:eta} and \eqref{eq:psi-h},
the reduced closure of
$\psi_1^{-1}(B_T\cap\{U\ne0\})$ is the closure
of the image of $\widehat C$ under the projection
$(u,s,\kappa)\mapsto(u,s)$.
The construction in Theorem~A with $h=1$
therefore recovers $C$ from $B_T$.

\section{Quotient surfaces and self-dual plane curves}

\subsection{Contact automorphisms of the incidence threefold}
Let $A:V\to V^\vee$ be an invertible linear map.
It defines a contact automorphism
\[
T_A:\mathcal I\longrightarrow\mathcal I,\qquad
(x,\ell)\longmapsto
(A_*^{-1}(\ell),(A^\vee)_*(x)).
\]
Indeed, $T_A$ is the composition of the exchange
of the two factors, which preserves the contact
structure, with the conormal lift of $A_*^{-1}$.

If $C\subset\Pp(V)$ is a reflexive plane curve
such that $C^\vee=A_*(C)$, then
\[
T_A(\widehat C)=\widehat C.
\]
Set
\[
G=A^{-1}A^\vee.
\]
Then
\[
T_A^2=\widehat G.
\]
In particular, $G_*(C)=C$.

\begin{lemma}\label{lem:normal-form}
Let $C\subset\Pp(V)$ be a reflexive self-dual plane
curve. There exist a basis of $V$ and an invertible
linear map $A:V\to V^\vee$ such that
$C^\vee=A_*(C)$ and the matrix of $A$ in this basis
and its dual basis is $A_h$, given by \eqref{eq:A-h},
for some $h\in k^\times$.
\end{lemma}

\begin{proof}
Choose an invertible linear map $A:V\to V^\vee$
such that $C^\vee=A_*(C)$ and put $G=A^{-1}A^\vee$.
Since $\det G=1$ and $A^\vee-A$ is represented
by a skew-symmetric matrix of odd size, the
eigenvalues of $G$ are $1,\lambda,\lambda^{-1}$
for some $\lambda\in k^\times$.
Moreover,
\[
G^\vee A G=A.
\]

We first show that $A$ can be chosen so that $G$
is semisimple. If $S\in\operatorname{GL}(V)$
satisfies
\[
S_*(C)=C,\qquad S^\vee A S=A,
\]
then replacing $A$ by $AS^{-1}$ preserves the
equality $C^\vee=A_*(C)$ and replaces $G$ by $S^2G$.

Suppose that $\operatorname{char} k=p>2$.
For sufficiently large $r$, the matrix $G^{p^r}$
is semisimple. Taking
\[
S=G^{(p^r-1)/2}
\]
gives $S^2G=G^{p^r}$.

Suppose now that $\operatorname{char} k=0$ and
that $G$ is not semisimple. Then $\lambda=\pm1$,
so $G^2$ is unipotent. Since $\dim V=3$, we have
$(G^2-I)^3=0$. Put
\[
S(t)=I+t(G^2-I)+\frac{t(t-1)}{2}(G^2-I)^2,
\qquad t\in k.
\]
For every integer $n$ we have $S(n)=(G^2)^n$.
Since
\[
S(n)_*(C)=C,\qquad S(n)^\vee A S(n)=A
\]
for every integer $n$, the same equalities hold
for all $t\in k$.
Taking $S=S(-1/4)$, we obtain
\[
(S^2G)^2=S^4G^2=I.
\]
Thus $S^2G$ is semisimple.

We may therefore assume that $G$ is semisimple.
If $G=I$, then $A$ is symmetric, and a suitable
choice of basis gives the matrix $A_1$.

Otherwise, choose $h\in k^\times$ such that
$h^2=\lambda$. The relation $G^\vee AG=A$ implies
that the eigenspaces of $G$ with eigenvalues
$\alpha$ and $\beta$ are orthogonal with respect
to the bilinear form represented by $A$ unless
$\alpha\beta=1$.

If $\lambda\ne-1$, choose a basis of eigenvectors
of $G$ with eigenvalues $h^{-2},1,h^2$,
respectively. If $\lambda=-1$, the relation
$A^\vee=AG$ implies that the restriction of the
bilinear form represented by $A$ to the
$(-1)$-eigenspace is alternating. In either case,
a suitable basis gives
\[
A=\begin{pmatrix}
0&0&a\\
0&b&0\\
c&0&0
\end{pmatrix},
\qquad
\frac ac=h^{-2},
\qquad abc\ne0.
\]
Rescaling the basis vectors, we obtain
$c=1$ and $b=-h^{-1}$, and hence $A=A_h$.
\end{proof}

In the homogeneous coordinates given by
Lemma~\ref{lem:normal-form}, consider the affine
chart \eqref{eq:incidence-chart} of $\mathcal I$ and
the change of coordinates \eqref{eq:affine-coordinates}.
These formulas are valid over $k$ since its
characteristic is different from $2$.
In the coordinates $(x,z,v)$,
\[
T_{A_h}(x,z,v)=\bigl(hv,h^2(xv-z),hx\bigr).
\]
In the coordinates $(u,s,\kappa)$, the transformation
$T_{A_h}$ takes the form
\begin{equation}\label{eq:diag}
T_{A_h}(u,s,\kappa)=(hu,-h^2s,-h\kappa).
\end{equation}

Let $C$ be a reflexive plane curve such that
$C^\vee=(A_h)_*(C)$. Since neither $C$ nor
$C^\vee$ is a line, the conormal curve
$\widehat C$ meets the affine chart above.
We regard $x,z,v,u,s,\kappa$ as rational functions
on $\widehat C$. Since the first projection
$\widehat C\to C$ is birational, we have
$k(\widehat C)=k(x,z)$. Moreover,
\[
\dd z=v\,\dd x.
\]
As $k$ is algebraically closed, the differentials
$\dd x$ and $\dd z$ cannot both vanish.
Hence $\dd x\ne0$.

Since $x=u+\kappa$, at least one of $\dd u$ and
$\dd\kappa$ is nonzero. The projective change
of coordinates
\[
[X:Y:Z]\longmapsto[X:Y:-Z]
\]
replaces the matrix $A_h$ by $A_{-h}$ up to
a nonzero scalar and induces the transformation
\[
(u,s,\kappa)\longmapsto
(\kappa,-s-2u\kappa,u).
\]
Thus, after a suitable choice of the basis and
the sign of $h$ in Lemma~\ref{lem:normal-form},
we may assume that $\dd u\ne0$ on $\widehat C$.

\subsection{The quotient surfaces \texorpdfstring{$Y_h$}{Yh}}
Let $h\in k^\times$ be a root of unity of order $m$.
The cyclic group $G_h=\langle R_h\rangle$ has order
$\operatorname{lcm}(m,2)$, which is invertible in $k$.
Consider the quotient surface
\[
Y_h=\A^2_{u,s}/G_h=\Spec k[u,s]^{G_h},
\]
and denote the quotient morphism by $\nu_h:\A^2_{u,s}\to Y_h$.
Since $G_h$ acts diagonally on $\A^2_{u,s}$,
the invariant ring is generated by monomials.
Thus $Y_h$ is an affine toric surface with at
most one cyclic quotient singularity.

On the open set $u\ne0$, put
\[
\xi=\frac{s}{u^2}.
\]
Then
\[
R_h(u,\xi)=(hu,-\xi).
\]
If $m$ is odd, the group $G_h$ is generated by
the transformations $(u,\xi)\mapsto(hu,\xi)$ and
$(u,\xi)\mapsto(u,-\xi)$. If $m$ is even, the
function $\xi/u^{m/2}$ is invariant under $G_h$.
In either case,
\begin{equation}\label{eq:quotcoord}
k[u,u^{-1},\xi]^{G_h}=k[U,U^{-1},V],
\end{equation}
where $U,V$ are given by \eqref{eq:psi-h}.

\begin{proposition}\label{prop:Yh}
The quotient surface $Y_h$ is rational.
The open subset $\{U\ne0\}\subset Y_h$ is
isomorphic to $\mathbb G_m\times\A^1$ by \eqref{eq:quotcoord}.
Moreover, $Y_1\simeq\A^2$, and
$Y_h\simeq\A^2$ if $h^2=-1$, whereas
\[
Y_{-1}\simeq\Spec k[A,B,C]/(AC-B^2).
\]
\end{proposition}

\begin{proof}
The first two assertions follow from the
description of the invariant ring on $u\ne0$.
For $h=1$, the invariant ring is $k[u,s^2]$.
If $h^2=-1$, it is $k[u^4,s]$.
For $h=-1$, the invariant ring is generated
by $A=u^2$, $B=us$, and $C=s^2$, with the
single relation $AC=B^2$.
\end{proof}

\subsection{Construction of reflexive self-dual plane curves}

Let $C$ be a reflexive self-dual plane curve.
By Lemma~\ref{lem:normal-form} and the choice
of coordinates in Section~2.1, we may assume that
$C^\vee=(A_h)_*(C)$ for some $h\in k^\times$
and that $\dd u\ne0$ on $\widehat C$.

By \eqref{eq:affine-contact-form}, the Legendrian
condition on $\widehat C$ gives $\dd s+2\kappa\,\dd u=0$,
and hence
\[
\kappa=-\frac12\frac{\dd s}{\dd u}.
\]
Since $\dd u\ne0$, the extension
$k(\widehat C)/k(u)$ is finite and separable.
The derivation $\dd/\dd u$ therefore restricts
to the subfield $k(u,s)$, so that
$\kappa\in k(u,s)$. The identities
\[
x=u+\kappa,\qquad
z=s+u\kappa+\frac{u^2-\kappa^2}{2}
\]
imply that
\[
k(\widehat C)=k(x,z)=k(u,s).
\]
Thus the rational map
\[
(u,s):\widehat C\dashrightarrow\A^2_{u,s}
\]
is birational onto the closure $D$ of its image.
By \eqref{eq:diag}, we have $R_h(D)=D$.

We now show that $h$ can be chosen to be a root
of unity. Suppose that $h$ is not a root of unity
and put $\xi=s/u^2\in k(D)$. The automorphism of
$k(D)$ induced by $R_h^2$ fixes $\xi$ and sends
$u$ to $h^2u$. If $\xi$ were nonconstant, this
would give an automorphism of infinite order
of the finite extension $k(D)/k(\xi)$, which is
impossible. Hence $\xi\in k$. Since
$\xi\circ R_h=-\xi$, we obtain $\xi=0$.
Consequently, $D=\{s=0\}$, $\kappa=0$, and
$C$ is the conic parametrized by
\[
[1:u:u^2/2].
\]
In this case $R_1(D)=D$, so we may replace
$h$ by $1$.

Conversely, let $h\in k^\times$ be a root of unity
and let $D\subset\A^2_{u,s}$ be an irreducible
curve such that $R_h(D)=D$ and $\dd u\ne0$
in $\Omega_{k(D)/k}$.
Let $X$ be the normalization of the projective
closure of $D$, and let $\sigma:X\to X$ be the
automorphism induced by the restriction of $R_h$
to $D$. Put
\[
\kappa=-\frac12\frac{\dd s}{\dd u}\in k(X).
\]
Consider the morphism $f:X\to\Pp^2$ given by
\begin{equation}\label{eq:param}
f=
\bigl[1:u+\kappa:s+u\kappa+\tfrac12(u^2-\kappa^2)\bigr].
\end{equation}
Write $f=[1:x:z]$ and put $v=u-\kappa$.
Then $\dd z=v\,\dd x$. Since
\[
u\circ\sigma=hu,\qquad
s\circ\sigma=-h^2s,
\]
we have $\kappa\circ\sigma=-h\kappa$, and hence
\[
x\circ\sigma=hv,\qquad v\circ\sigma=hx.
\]
If $\dd x=0$, then $\dd v=0$ and consequently
$2\dd u=\dd x+\dd v=0$, a contradiction.
Thus $\dd x\ne0$ and $\dd v\ne0$.

Let $C=f(X)$. Since $\dd x\ne0$, the extension
$k(X)/k(x)$ is finite and separable.
The derivation $\dd/\dd x$ therefore restricts
to $k(x,z)$, so that
\[
v=\frac{\dd z}{\dd x}\in k(x,z).
\]
The formulas \eqref{eq:affine-coordinates} give
$k(X)=k(u,s)=k(x,z)=k(C)$.
Hence $f$ is birational onto $C$.

Since $v$ is nonconstant, $C$ is not a line.
The morphism $f^\vee:X\to\Pp(V^\vee)$ obtained
by composing $f$ with the Gauss map of $C$
on the inverse image of its smooth locus is
given by
\[
f^\vee=[xv-z:-v:1].
\]
A direct computation gives
\[
f^\vee=(A_h)_*\circ f\circ\sigma.
\]
Consequently, $f^\vee$ is birational onto its
image and $C^\vee=(A_h)_*(C)$.
Thus $C$ is reflexive and self-dual.

\begin{proposition}\label{prop:quotientA}
Let $h\in k^\times$ be a root of unity, and let
$D\subset\A^2_{u,s}$ be an irreducible closed curve.
Then $R_h(D)=D$ if and only if
\[
D=\nu_h^{-1}(B)_{\mathrm{red}}
\]
for some irreducible closed curve $B\subset Y_h$.
\end{proposition}

\begin{proof}
The reduced inverse image of a curve in $Y_h$ is
$R_h$-invariant. Conversely, if $R_h(D)=D$, put
$B=\nu_h(D)$. Since $\nu_h$ is finite and its
fibres over closed points are precisely the
$G_h$-orbits, $B$ is an irreducible closed curve
and $D=\nu_h^{-1}(B)_{\mathrm{red}}$.
\end{proof}

Since $\dd u\ne0$, each curve $D$ considered above
meets the open set $u\ne0$. On this open set,
the quotient morphism $\nu_h$ is given by
\eqref{eq:psi-h}. By Proposition~\ref{prop:quotientA},
the intersection of an $R_h$-invariant curve $D$
with this open set is the reduced inverse image
of an irreducible closed curve in $\{U\ne0\}$.
Taking its closure in $\A^2_{U,V}$ gives the curve
$B$ in Theorem~A. Conversely, the reduced closure
of $\psi_h^{-1}(B\cap\{U\ne0\})$ is $R_h$-invariant.
Theorem~A now follows from the construction
\eqref{eq:param} and its converse.

\begin{remark}\label{rem:irred}
Both the irreducibility condition and the condition
$\dd u\ne0$ in Theorem~A can be tested on the
function field $k(B)$.

If $m$ is even, the reduced inverse image $D$
is irreducible if and only if $U$ is not an
$r$-th power in $k(B)$ for any prime divisor
$r$ of $m$.

If $m$ is odd, the same condition on $U$ is
necessary, and, provided that $V\ne0$ in $k(B)$,
one must also require that $V$ is not a square
in $k(B)$. These conditions are sufficient
because $m$ is odd.

Indeed, in the coordinates $(u,\xi)$ the inverse
image is defined over $k(B)$ by
\[
\begin{cases}
u^m=U,\quad \xi^2=V,&\text{if $m$ is odd},\\
u^m=U,\quad \xi=Vu^{m/2},&\text{if $m$ is even}.
\end{cases}
\]
The conditions on $U$ and $V$ for $D$ to be
irreducible follow from Kummer theory.
If $m$ is odd and $V=0$, the reduced
inverse image is defined by $u^m=U$ and $\xi=0$.

Assuming that $D$ is irreducible, the extension
$k(D)/k(B)$ is finite and separable.
Since $U=u^m$ and $m$ is invertible in $k$, the
condition $\dd u\ne0$ in $\Omega_{k(D)/k}$ is
equivalent to $\dd U\ne0$ in $\Omega_{k(B)/k}$.
\end{remark}

\subsection{A curve not self-dual under any polarity}

\begin{proposition}\label{prop:nonpolar}
There exists an irreducible reflexive self-dual
plane curve whose normalization has genus $15$
and which is not self-dual under any polarity.
\end{proposition}

\begin{proof}
Choose $h\in k$ such that $h^2=-1$.
Let $P(s)\in k[s]$ be a polynomial of degree $12$
with simple roots such that no nonidentity
automorphism of $\Pp^1$ preserves the set of
its roots. Such a polynomial exists: for each
nonidentity permutation of twelve distinct points
of $\Pp^1$, there is a set of four points which
is not preserved by the permutation. The equality
of the cross-ratios of the corresponding ordered
quadruples imposes a nontrivial equation on the
configuration of the twelve points.

Let $D\subset\A^2_{u,s}$ be the irreducible curve
defined by
\[
u^4=P(s),
\]
and let $X$ be the normalization of its projective
closure. The degree-four morphism $s:X\to\Pp^1$
is totally ramified over the twelve roots of $P$
and unramified at infinity. By the Riemann--Hurwitz
formula,
\[
2g(X)-2=-8+12\cdot3,
\]
so $g(X)=15$.

Since $h^2=-1$, we have $R_h(u,s)=(hu,s)$.
Let $\eta:X\to X$ be the automorphism induced
by the restriction of $R_h$ to $D$.
The extension $k(X)/k(s)$ is separable, so
$\dd s\ne0$. Differentiating the equation of $D$,
we obtain
\[
4u^3\dd u=P'(s)\dd s.
\]
Since $P$ has simple roots, $\dd u\ne0$.
By Theorem~A, the morphism $f:X\to\Pp^2$
given by \eqref{eq:param} is birational onto
a reflexive self-dual plane curve $C$.

We show that every degree-four morphism
$t:X\to\Pp^1$ differs from $s$ by an automorphism
of $\Pp^1$. Put
\[
E=k(s,t)\subset k(X),\qquad e=[k(X):E].
\]
Since both $s$ and $t$ have degree four,
$e\in\{1,2,4\}$. The image of the morphism
\[
(s,t):X\longrightarrow\Pp^1\times\Pp^1
\]
is a curve of bidegree $(4/e,4/e)$ with function
field $E$. Hence its normalization has genus at
most $(4/e-1)^2$.

If $e=1$, this gives $g(X)\le9$, a contradiction.
If $e=2$, then $E=k(s,u^2)$, since $k(X)/k(s)$
is a cyclic extension of degree four.
The smooth projective curve with function field
$k(s,u^2)$ is a double cover of $\Pp^1$ branched
over the twelve roots of $P$, and thus has genus
$5$. On the other hand, the genus bound for
$e=2$ gives $g\le1$, a contradiction.
Consequently, $e=4$ and $k(s)=k(t)$.

It follows that every automorphism of $X$
descends via $s$ to an automorphism of $\Pp^1$
preserving the set of roots of $P$. By the
choice of $P$, this automorphism is the identity.
Therefore
\[
\Aut(X)=\langle\eta\rangle
\simeq\mathbb Z/4\mathbb Z.
\]

Suppose that $C$ is self-dual under a polarity
represented by an invertible symmetric matrix $Q$.
The contact automorphism $T_Q$ induces an
automorphism $\tau$ of $X$ satisfying
$\tau^2=\id_X$ and
\[
f^\vee=Q_*\circ f\circ\tau.
\]
Thus $\tau=\id_X$ or $\tau=\eta^2$.

By \eqref{eq:param},
\[
f\circ\eta^2=S_*\circ f,\qquad
S=\begin{pmatrix}
1&0&0\\
0&-1&0\\
0&0&1
\end{pmatrix}.
\]
It follows that
\[
f^\vee=(QS^j)_*\circ f
\]
for some $j\in\{0,1\}$. Since each smooth
point of $C$ belongs to its tangent line,
\[
p^TQS^jp=0
\]
for every nonzero representative $p\in k^3$
of a point of $C$.

The quadratic form $p\mapsto p^TQS^jp$ is
nonzero, since an invertible matrix of odd size
cannot be skew-symmetric.
Hence $C$ is contained in a conic, contradicting
$g(X)=15$.
\end{proof}

\subsection{Self-dual curves of large genus}

\begin{proof}[Proof of Proposition~\ref{prop:high-genus}]
Fix a positive integer $b$ divisible by $8$ and put
\[
\ell=4b+1,\qquad
N=\frac b4\ell(\ell-2),\qquad
a=\frac N\ell=b^2-\frac b4.
\]
For $1\le j\le 2b$, put
\[
\varepsilon_j=
\begin{cases}
1,&1\le j\le b,\\
-1,&b<j\le 2b.
\end{cases}
\]
In the affine coordinate $w$ on $\Pp^1$, consider
the points
\[
r_0=-\frac{\ell}{\ell-2},\qquad
r_j=r_0\frac{2j}{b+2j},\quad 1\le j\le 2b,
\]
and the rational function
\[
F(w)=\prod_{j=1}^{2b}(w-r_j)^{\varepsilon_j}.
\]
The numbers $r_0,\ldots,r_{2b}$ are distinct
and nonzero. Put
\[
e_0=-\frac b4,\qquad
e_j=\varepsilon_j(j-a),\quad 1\le j\le 2b.
\]
Since $a>2b$, all the integers $e_j$ are nonzero.

Let $\kappa$ be a root of the equation
\begin{equation}\label{eq:high-genus-cover}
\kappa^N=\prod_{j=0}^{2b}(w-r_j)^{e_j}
\end{equation}
in an algebraic closure of $\C(w)$, and let $X_b$
be the smooth projective curve with function field
$\C(w,\kappa)$.

The extension $\C(w,\kappa)/\C(w)$ is cyclic,
since $\C$ contains all $N$-th roots of unity.
If its degree is $m$, then $m$ divides $N$
and $\kappa^m\in\C(w)$. Taking the orders at
$r_j$ in the equality
\[
(\kappa^m)^{N/m}
=\prod_{j=0}^{2b}(w-r_j)^{e_j},
\]
we see that $N/m$ divides every $e_j$.
Since $e_1$ and $e_2$ are consecutive integers,
we obtain $m=N$.

Consider the rational functions on $X_b$ given by
\[
t=F(w)\kappa^\ell,\qquad
u=\frac{tw}{\kappa},\qquad
s=t-u\kappa.
\]
We claim that
\begin{equation}\label{eq:high-genus-differential}
\dd t=u\,\dd\kappa-\kappa\,\dd u.
\end{equation}
Substituting the expressions for $t$ and $u$,
we see that \eqref{eq:high-genus-differential}
is equivalent to
\begin{equation}\label{eq:high-genus-logarithmic}
\frac{\dd\kappa}{\kappa}
=-\frac{1+w}{\ell+(\ell-2)w}\frac{\dd F}{F}
-\frac{\dd w}{\ell+(\ell-2)w}.
\end{equation}

By \eqref{eq:high-genus-cover},
\[
\frac{\dd\kappa}{\kappa}
=\frac1N\sum_{j=0}^{2b}
e_j\frac{\dd w}{w-r_j}.
\]
Thus both sides of \eqref{eq:high-genus-logarithmic}
are pullbacks of rational differential forms on
$\Pp^1$ with at most simple poles at
$r_0,\ldots,r_{2b},\infty$.

For $1\le j\le 2b$, we have
\[
-\frac{1+r_j}{\ell+(\ell-2)r_j}
=-\frac1\ell+\frac jN.
\]
Hence the residues of both differential forms
at $r_j$ are equal to
\[
\varepsilon_j\left(-\frac1\ell+\frac jN\right)
=\frac{e_j}{N}.
\]
Since $F$ is regular at infinity and $F(\infty)=1$,
the residue of the differential form on the
right-hand side of \eqref{eq:high-genus-logarithmic}
at infinity is $1/(\ell-2)$. Put
\[
e_\infty=-\sum_{j=0}^{2b}e_j
=b^2+\frac b4=\frac N{\ell-2}.
\]
The residue of the differential form on the
left-hand side at infinity is therefore also
$1/(\ell-2)$.
By the residue theorem, the residues at $r_0$
coincide as well. The difference of the two
differential forms is regular on $\Pp^1$ and
hence vanishes. This proves
\eqref{eq:high-genus-logarithmic} and
\eqref{eq:high-genus-differential}.

Since $s=t-u\kappa$, we obtain
\[
\dd s=-2\kappa\,\dd u.
\]
The function $u=wF(w)\kappa^{\ell-1}$ has a pole
at every point of $X_b$ over $r_0$, since
$e_0<0$ and both $w$ and $F(w)$ are regular
and nonvanishing at $r_0$. In particular,
$\dd u\ne0$, so that
\[
\kappa=-\frac12\frac{\dd s}{\dd u}
\in\C(u,s).
\]
Since $t=s+u\kappa$ and $w=u\kappa/t$, we have
\[
\C(X_b)=\C(u,s).
\]

The integer $N$ is even and $\ell$ is odd.
Consequently, the automorphism $\sigma$ of $X_b$
given by
\[
w\circ\sigma=w,\qquad
\kappa\circ\sigma=-\kappa
\]
is an involution satisfying
\[
u\circ\sigma=u,\qquad s\circ\sigma=-s.
\]
Let $D\subset\A^2_{u,s}$ be the closure of the
image of the rational map $(u,s)$.
Then $R_1(D)=D$, and $X_b$ is isomorphic to
the normalization of the projective closure
of $D$.

By Theorem~A, the morphism $f:X_b\to\Pp^2$
given by \eqref{eq:param} is birational onto
a plane curve $C_b$ satisfying
\[
C_b^\vee=(A_1)_*(C_b).
\]
Explicitly, we have
\[
f=[1:x:z],\qquad
x=u+\kappa,\qquad
z=t+\frac{u^2-\kappa^2}{2}.
\]

Put
\[
h_j=\gcd(N,|e_j|),\qquad 0\le j\le 2b,
\qquad
h_\infty=\gcd(N,e_\infty)=e_\infty.
\]
The morphism $w:X_b\to\Pp^1$ has $h_j$ points
over $r_j$, each of ramification index $N/h_j$.
At each such point $P$, we have
\[
\operatorname{ord}_P(w-r_j)=\frac{N}{h_j},
\qquad
\operatorname{ord}_P(\kappa)=\frac{e_j}{h_j}.
\]
Indeed, after absorbing a local unit into the
$N$-th power, equation \eqref{eq:high-genus-cover}
takes the form $\kappa^N=(w-r_j)^{e_j}$.
Similarly, there are $h_\infty$ points over infinity,
each of ramification index $N/h_\infty$.
The morphism $w$ is unramified elsewhere.

Using the expressions for $t$ and $u$, we obtain
the following orders at a point $P$ over $r_j$:
\[
\begin{array}{c|ccc}
 & h_j\operatorname{ord}_P(\kappa)
 & h_j\operatorname{ord}_P(t)
 & h_j\operatorname{ord}_P(u)\\ \hline
1\le j\le b
 & j-a & \ell j & a+(\ell-1)j\\
b<j\le 2b
 & a-j & -\ell j & -a-(\ell-1)j\\
j=0
 & -b/4 & -\ell b/4 & -(\ell-1)b/4
\end{array}
\]
At each point over infinity, the orders of
$\kappa,t,u$ are $1,\ell,1$, respectively, since
$w$ has a pole of order $\ell-2$ and $F(\infty)=1$.
The functions $\kappa,t,u$ are regular at all
other points of $X_b$.

For $1\le j\le b$, only $\kappa$ has a pole
at the points over $r_j$. Hence the pole orders
of $x$ and $z$ at each such point are
\[
\frac{a-j}{h_j},\qquad
\frac{2(a-j)}{h_j},
\]
respectively.

For $b<j\le 2b$, the term $u^2/2$ has a pole
of strictly larger order than the term $t$,
since
\[
2(a+(\ell-1)j)>\ell j.
\]
Thus the pole orders of $x$ and $z$ at each
point over $r_j$ are
\[
\frac{a+(\ell-1)j}{h_j},\qquad
\frac{2(a+(\ell-1)j)}{h_j},
\]
respectively.

At each point over $r_0$, the term $u^2/2$
has a pole of strictly larger order than
either $t$ or $\kappa^2/2$. Hence the pole
orders of $x$ and $z$ are
\[
\frac{(\ell-1)b}{4h_0},\qquad
\frac{(\ell-1)b}{2h_0},
\]
respectively.

Let $E$ be the pole divisor of $z$.
Since the pole divisor of $x$ is bounded by $E$,
the functions $1,x,z$ determine sections of
$\mathcal O_{X_b}(E)$ with no common zero.
Consequently,
\[
f^*\mathcal O_{\Pp^2}(1)\simeq
\mathcal O_{X_b}(E).
\]
Since $f$ is birational onto $C_b$, we obtain
\begin{align*}
d_b=\deg E
&=2\left(
\sum_{j=1}^{b}(a-j)
+\sum_{j=b+1}^{2b}(a+(\ell-1)j)
+\frac{(\ell-1)b}{4}
\right)\\
&=16b^3+4b^2-b.
\end{align*}

By the Riemann--Hurwitz formula,
\begin{align*}
2g_b-2
&=-2N+\sum_{j=0}^{2b}(N-h_j)+(N-h_\infty)\\
&=2bN-\sum_{j=0}^{2b}h_j-h_\infty.
\end{align*}
Since $h_j\le |e_j|=O(b^2)$ and
$h_\infty=e_\infty=O(b^2)$, we have
\[
\sum_{j=0}^{2b}h_j+h_\infty=O(b^3).
\]
Using $N=4b^3-b/4$, we obtain
\[
g_b=4b^4+O(b^3).
\]
Since $d_b=16b^3+O(b^2)$, it follows that
\[
g_b=2^{-10/3}d_b^{4/3}+O(d_b).
\]
\end{proof}

\section{Normalizations of self-dual plane curves}

\begin{proof}[Proof of Theorem~B]
Assume that $X$ admits a nontrivial involution
$\sigma$. Put
\[
F=k(X),\qquad E=F^{\langle\sigma\rangle}.
\]
Since $F/E$ is a separable quadratic extension,
\[
\Omega_{F/k}\simeq
F\otimes_E\Omega_{E/k}.
\]
Thus there exists $u\in E$ such that $\dd u\ne0$.
In particular, $F/k(u)$ is finite and separable.

Consider the $k(u)$-vector subspace
\[
F^-=\{s\in F:s\circ\sigma=-s\}.
\]
For any nonzero $s_0\in F^-$ we have
$F^-=s_0E$ and $F=E(s_0)$.
Consequently, $F^-$ generates $F$ as a field
over $k(u)$. Any two distinct $k(u)$-embeddings
of $F$ into an algebraic closure of $k(u)$
therefore have distinct restrictions to $F^-$.
Since $k(u)$ is infinite, we can choose
$s\in F^-$ outside the finitely many proper
linear subspaces on which two such embeddings
coincide. Then $F=k(u,s)$.

We have thus obtained rational functions $u,s$
on $X$ satisfying
\[
k(X)=k(u,s),\qquad
u\circ\sigma=u,\qquad
s\circ\sigma=-s,\qquad
\dd u\ne0.
\]
Set
\[
\kappa=-\frac12\frac{\dd s}{\dd u}
\]
and define $f:X\to\Pp^2$ by \eqref{eq:param}.
Let $\bar k$ be an algebraic closure of $k$.
Since $X$ is geometrically integral, we have
$\bar k(X_{\bar k})=\bar k(u,s)$.
Moreover, $\dd u\ne0$ in
$\Omega_{\bar k(X_{\bar k})/\bar k}$,
and the identities $u\circ\sigma=u$ and
$s\circ\sigma=-s$ also hold on $X_{\bar k}$.
Applying the argument of Section~2.3 following
\eqref{eq:param} with $h=1$, and using the fact
that birationality and reflexivity can be checked
over $\bar k$, we conclude that $f$ is birational
onto a geometrically integral reflexive plane
curve $C$ of degree at least two.
The same argument gives the identity
\[
f^\vee=\rho_B\circ f\circ\sigma
\]
over $\bar k$ and hence over $k$, where $B$ is
the nondegenerate symmetric bilinear form
with $b_B=A_1$.
This proves \textup{(1)}\,$\Rightarrow$\,\textup{(3)}
and the last assertion of the theorem.
The implication
\textup{(3)}\,$\Rightarrow$\,\textup{(2)}
is immediate.

Now assume \textup{(2)} and let $C$ be a
geometrically integral reflexive self-dual plane
curve with normalization $X$. We may
assume that $C$ is not a conic.

Choose a linear isomorphism $A:V\to V^\vee$
such that $C^\vee=A_*(C)$.
The replacement of $A$ by $AS^{-1}$ in the proof
of Lemma~\ref{lem:normal-form} is defined over $k$.
Thus we may assume that $G=A^{-1}A^\vee$ is semisimple.

The automorphism $T_A$ of $\mathcal I$ preserves
$\widehat C$ and induces an automorphism
$\tau:X\to X$.
Let $\bar k$ be an algebraic closure of $k$.
Over $\bar k$, choose a basis of $V\otimes_k\bar k$
in which $A$ has matrix $A_h$ and $\dd u\ne0$,
as in Section~2.1.
By the first part of Section~2.3, we have
$\bar k(X)=\bar k(u,s)$ and, since $C$ is not
a conic, $h$ is a root of unity and $s\ne0$.
By \eqref{eq:diag},
\[
u\circ\tau=hu,\qquad
s\circ\tau=-h^2s.
\]

Let $m$ be the order of $h$.
If $m$ is odd, then $\tau^m$ is a nontrivial
involution, since
\[
u\circ\tau^m=u,\qquad
s\circ\tau^m=-s.
\]
If $m$ is even, then $\tau^{m/2}$ is a
nontrivial involution, since
\[
u\circ\tau^{m/2}=-u,\qquad
s\circ\tau^{m/2}=(-1)^{m/2}s.
\]
As $\tau$ is defined over $k$, so are its powers.
This proves \textup{(2)}\,$\Rightarrow$\,\textup{(1)}.
\end{proof}

\printbibliography

\end{document}